\input amstex

\NoBlackBoxes

\documentstyle{amsppt}

\topmatter

\title{Multidimensional topological Galois theory}\endtitle

\author Askold Khovanskii\endauthor
\affil{Department of Mathematics, University of Toronto, Toronto, Canada}\endaffil

\thanks{This work was partially supported by the Canadian Grant No. 156833-17.}\endthanks

\abstract{In this preprint we present an outline of the multidimensional  version of topological Galois theory. The theory   studies topological obstruction to solvability  of equations ``in finite terms" (i.e. to their solvability by radicals, by elementary functions, by quadratures and so on).
This preprint is based on the author's book on topological Galois theory. It contains definitions, statements of results and comments to them. Basically no proofs are presented.

 This preprint was  written as a part of the comments to a new edition (in preparation) of the classical book ``Integration in finite terms'' by J.F.~Ritt.  }\endabstract



\endtopmatter
\document

\subhead{1. Introduction}\endsubhead
In topological Galois theory for functions of one variable (see [1], [2]), it is proved that
the way the Riemann surface of a function is positioned over the complex line can obstruct
the representability of this function ``in finite term'' (i.e. its representability by radicals, by quadratures, by generalized quadratures and so on).
This not only explains why many algebraic and differential equations are not solvable in finite terms,
but also gives the strongest known results on their unsolvability.

In the multidimensional version of topological Galois theory  analogous results are proved. But in the multidimensional case all constructions and proofs are much more complicated and involved than in the one dimensional case (see [1]).

\subhead{2. Classes of functions}
\endsubhead
An equation is solvable ``in finite terms'' (or is solvable ``explicitly'') if its solutions belong to a certain class of functions. Different classes of functions correspond to different notions of  solvability in finite terms.

A class of functions can be introduced by specifying a list of {\it basic functions}
and a list of {\it admissible operations}.
Given the two lists, the class of functions is defined as the set
of all functions that can be obtained
from the basic functions by repeated application of admissible operations.
Below, we define Liouvillian classes of functions  in exactly this way.

Classes of functions, which appear in the problems of integrability in finite terms,
contain multivalued functions.
Thus the basic terminology should be made clear.

We understand operations on multivalued functions
of several variables  in a slightly more
restrictive sense than operations on multivalued functions of
single variable (one dimensional case is discussed  in [1], [2].

 Fix a class of basic functions and some set of admissible operations. Can a given
function (which is obtained, say, by solving a certain
algebraic or a differential equation) be expressed through the
basic functions by means of admissible operations? We are
interested in various {\it single valued branches} of
multivalued functions over various domains. Every function,
even if it is multivalued, will be considered as a collection
of all its single valued branches. We will only apply
admissible operations (such as arithmetic operations and
composition) to single valued branches of the function over
various domains. Since we deal with analytic functions, it
suffices only to consider small neighborhoods of points as
domains.

We can now rephrase the question in the following way:
{\it can a given function germ at a given point be expressed through the germs of basic functions
with the help of admissible operations?}
Of course, the answer depends on the choice of a point and on the choice of a single valued
germ at this point belonging to the given multivalued function.
It turns out, however, that for the classes of functions interesting to us the desired expression
is either impossible for every germ of a given multivalued function at every point or
the ``same'' expression serves all germs of a given multivalued function at almost every point
of the space.

For  functions of one variable, we use a different, extended definition of operations on
multivalued functions, in which the multivalued function is viewed as a single object. This definition is essentially equivalent to including the operation of analytic continuation in the list of admissible operations on analytic germs (all details can be found in [1]).
For functions of many variables, we need to adopt the more restrictive understanding
of operations on multivalued functions, which is, however, no less (and perhaps even more) natural.
\subhead{3. Specifics of the multidimensional case}\endsubhead
I was always under impression that a full-fledged
multidimensional version of topological Galois theory was
impossible. The reason was that, to construct such a version
for the case of many variables, one would need to have 
information on extendability of function germs not only outside
their ramification sets but also along these sets. It seemed
that there was nothing to extract such information from.

To illustrate the problem  consider the following  situation. Let $f$ be a multivalued analytic function on $\Bbb{C}^n$, whose set of singular points is an analytic set $\Sigma_f\subset \Bbb C^n$. Let $ f_a$ be an analytic  germ of  $f$ at a point $a\in \Bbb C^n$. Let $g:(\Bbb C^k,b)\to
(\Bbb{C}^n,a)$ be an analytic map. Consider a germ $\varphi_b$ at the point $b\in \Bbb C^k$ of the composition $ f_a \circ g_b $. One can ask the following questions:

1) Is it true that $\varphi_b$  is a germ of a
multivalued function $\varphi$ on $\Bbb C^k$, whose set of singular points $\Sigma_\varphi$ is contained in a proper analytic subset of $\Bbb C^k$?

2)  Is it true that  the monodromy group $M_\varphi$ of $\varphi$  corresponding to motions  around the  set $\Sigma _\varphi \subset \Bbb C^k$ can be estimated in terns of the monodromy group $M_f$ of $f$  corresponding to motions  around the  set $\Sigma_f\subset \Bbb C^n $? For example, if $M_f$ is a solvable group is it true that $M_\varphi$ also is a solvable group?

If the image $g(\Bbb C^k)$ is not contained in the singular set $\Sigma_f$  then the answers to the both questions are positive: the set $\Sigma _\varphi $ belongs to the analytic set $g^{-1}(\Sigma_f)$ and the group $M_\varphi$ is a subgroup of a certain factor group of $M_f$. These statements are not complicated and can be proved by the same arguments as in the one dimensional topological Galois theory.

Assume that the multivalued function $f$ has an analytic germ  $f_a$ at a point $a$ belonging to the singular set $\Sigma_f$ (some of the germs of the multivalued function $f$ may
appear to be nonsingular at singular points of this function). Assume now that the image $g(\Bbb C^k)$ is  contained in the singular set $\Sigma_f$ and $a=g(b)$. It turns out that for the germ $\varphi_b=f_a\circ g_b$ the answers to the both questions also are positive.
In  this situation all the proofs are  more involved. They use new arguments from multidimensional complex analysis  and from  group theory.

It turns out that  function germs can
sometimes be automatically extended along their ramification
sets (see [1]). That new statement from complex analysis suggests the positive answer to the first question.

To describe the connection between the monodromy group of the
function $f$ and the monodromy groups of the composition $\varphi = f\circ g$, we
introduce and develop the notion of pullback closure for groups (see [1]). The use of this operation, in turn, forces us
to reconsider all arguments we used in the one dimensional version of topological Galois theory. As a result  we obtain a positive answer to the second question.

\bigskip
\subhead {4. Liouvillian classes of multivariate functions}\endsubhead
In this section we define Liouvillian classes of functions  for the case of several variables.
These classes  are defined in the same way as the corresponding classes for functions of one variable (see [1], [2]).
The only difference is in the details.

We fix an ascending chain of standard coordinate subspaces of
strictly increasing dimension:
$0\subset\Bbb C^1\subset\dots\subset\Bbb C^n\subset\dots$
with coordinate functions $x_1$, $\dots$, $x_n$, $\dots$ (for
every $k>0$, the functions $x_1$, $\dots$, $x_k$ are coordinate
functions on $\Bbb C^k$).
Below, we define Liouvillian
classes of functions for each of the standard coordinate
subspaces $\Bbb C^k$.

To define Liouvillian classes, we will need the list of basic elementary functions and the list of classical operations.
\medskip

{\bf List of basic elementary functions.}

1. All complex constants and all coordinate functions
    $x_1$, $\dots$, $x_n$ for every standard coordinate
    subspace $\Bbb C^n$.

2. The exponential, the logarithm and the power
    $x^\alpha$, where $\alpha$ is any complex constant.

3. Trigonometric functions:  sine, cosine, tangent,
    cotangent.

4. Inverse trigonometric functions: arcsine, arccosine,
    arctangent, arccotangent.
    \medskip

Let us now turn to the list of classical operations on
functions.

\medskip

{\bf List of classical operations.}

1. {\it Operation of composition} that takes a function $f$
    of $k$ variables and functions $g_1$, $\dots$, $g_k$ of $n$
    variables to the function $f(g_1,\dots,g_k)$ of $n$
    variables.

2. {\it Arithmetic operations} that take functions $f$
    and $g$ to the functions $f+g$, $f-g$, $fg$ and $f/g$.

3. {\it Operations of partial differentiation with
    respect to independent variables}. For functions of $n$
    variables, there are $n$ such operations: the $i$-th
    operation assigns the function $\frac{\partial
    f}{\partial x_i}$ to a function $f$ of the variables
    $x_1$, $\dots$, $x_n$.

4. {\it Operation of integration} that takes $k$
    functions $f_1$, $\dots$, $f_k$ of the variables $x_1$,
    $\dots$, $x_n$, for which the differential one-form
    $\alpha=f_1 dx_1+\dots+f_k dx_k$ is closed, to the
    indefinite integral $y$ of the form $\alpha$ (i.e. to
    any function $y$ such that $dy=\alpha$). The function
    $y$ is determined by the functions $f_1$, $\dots$,
    $f_k$ up to an additive constant.

5. {\it Operation of solving an algebraic equation} that
    takes functions $f_1,\dots,f_n$ to the function $y$
    such that $y^n+f_1y^{n-1}+\dots+f_n=0$. The function
    $y$ may not be quite uniquely determined by the
    functions $f_1$, $\dots$, $f_n$, since an algebraic
    equation of degree $n$ can have $n$ solutions.
\smallskip

We now resume defining Liouvillian classes of functions.
\medskip
{\bf Functions of $n$ variables representable by radicals.}
List of basic functions: All complex constants and all coordinate functions.
List of admissible operations:  composition,
arithmetic operations and the operation of taking the $m$-th root
$f^{\frac {1}{m}}$, $m=2,3,\dots$, of a given function~$f$.
\medskip

\medskip

{\bf Functions of $n$ variables representable by
$k$-radicals.} This class of functions is defined in the same
way as the class of functions representable by radicals. We
only need to add the operation of solving algebraic equations
of degree $\leq k$ to the list of admissible operations.

\medskip

{\bf Elementary functions of $n$ variables.}
List of basic functions: basic elementary functions.
List of admissible operations: composition, arithmetic
operations, differentiation.
\medskip
\medskip

{\bf Generalized elementary functions of $n$
variables.} This class of functions is defined in the same way
as the class of elementary functions. We only need to add the
operation of solving algebraic equations to the list of
admissible operations.

\medskip

{\bf Functions of $n$ variables representable by
quadratures.}
List of basic functions: basic elementary functions.
List of admissible operations: composition, arithmetic
operations, differentiation, integration.
\medskip

{\bf Functions of $n$ variables representable by
$k$-quadratures.} This class of functions is defined in the
same way as the class of functions representable by
quadratures. We only need to add the operation of solving
algebraic equations of degree at most $k$ to the list of
admissible operations.

\medskip
{\bf Functions $n$ variables representable by
generalized quadratures.} This class of functions is defined in
the same way as the class of functions representable by
quadratures. We only need to add the operation of solving
algebraic equations to the list of admissible operations.

\subhead {5. Strong non representability in finite terms}\endsubhead
Topological obstructions to  the representability of functions in finite terms relate to
branching. It turns out that if a function does not belong to a certain Liouvillian class by topological reasons then it automatically does not belong to a much wider {\it extended Liouvillian class of functions}.

Such an extended Liouvillian class is defined as follows: its list of admissible operations is the same as in the original Liouvillian class and its list of basic functions is the  list of basic function in the original class extended by all single valued functions of any number of variables having proper analytic set of singular points.

\definition {Definition} A germ $f$ is a germ of function belonging to the {\it extended class of functions representable by  by quadratures}  if it can be represented  by germs of basic elementary functions and by germs of single valued functions, whose set of singular points is a proper analytic set,  by means of composition, integration,  arithmetic  operations and
differentiation.

\enddefinition
\definition {Definition} A germ $f$ is {\it strongly non representable by quadratures}  if it is not a germ of function from the extended class of functions representable by  by quadratures.
\enddefinition

The definition of {\it strong non representability of a germ $f$} by radicals, by $k$-radical, by elementary functions, by generalized elementary functions, by $k$-quadratures and by generalized quadratures is similar to the above definition.

\subhead{6. Holonomic systems of linear differential equations}\endsubhead
Consider a system of $N$ linear
differential equations $L_j(y)=0$, $j=1,\dots, N$,
$$
L_j(y)=\sum a_{i_1,\dots, i_n}\frac{\partial^{i_1+\dots
+i_n}y}{\partial x_1^{i_1} \dots
\partial x^{i_n}_n} =0, \tag 1
$$
on an unknown function $y$, whose coefficients $a_{i_1,\dots, i_n}$
are analytic  functions in a domain $U\subset \Bbb C^n$.

The system (1) is {\it holonomic} if at every point $a\in U$ the $\Bbb C$-linear space  $V_a$ of germs $y_a$ satisfying the system (1) has finite dimension, $\dim_{\Bbb C}V_a=d(a)<\infty$.
Holonomic systems  can be considered as a multidimensional generalization of linear differential equation  on one unknown function of a single variable. Kolchin obtained a generalization of the Picard--Vessiot theory (Galois theory for linear differential equations)
to the case of holonomic systems of differential equations [3].

Holonomic system (1) has the following properties:

1) There exists an analytic {\it singular hypersurface  $\Sigma \subset U$} such that the dimension $d(a)=\dim_{\Bbb C}V_a$ is constant $d(a) \equiv d$ on $U\setminus \Sigma$.

2) Let $\gamma:I\to U\setminus \Sigma$ be a continuous map,  where $I$ is
 the unit segment $0\leq t\leq 1$ and $\gamma(0)=a$, $\gamma(1)=b$. Then the space $V_a$ of solutions of (1) at the point $a$  admits  analytic continuation  along $\gamma$ and the space obtained by the continuation at the point $b$ is the space $V_b$ of solutions of (1)  at the point $b$.

3) If all  equations of the system (1) admit analytic continuation  to some domain $W$, then the system  obtained by such a continuation  is a holonomic system  in the domain $W$.

Let $a\notin \Sigma$ be a point not belonging to the hypersurface $\Sigma$.
Take an arbitrary path $\gamma(t)$ in the domain
$U$ originating and terminating at $a$ and avoiding the
hypersurface $\Sigma$.
Solutions of this system admit analytic continuations along the path $\gamma$,
which are also solutions of the system.
Therefore, every such path $\gamma$ gives rise to a linear map $M_{\gamma}$
of the solution space $V_a$ to itself.
The collection of linear transformations $M_{\gamma}$ corresponding to all paths
$\gamma$ form a group, which is called the
{\it monodromy group of the holonomic system}.

\subhead{7. $\Cal {SC}$-germs}\endsubhead
There is a wide class of $\Cal S$-functions in one variable  containing all Liouvillian functions and stable under classical operations, for which the
monodromy group is defined. The class of $\Cal S$-functions plays an important role in the one dimensional  version of topological Galois theory (see [1], [2]). Is there a sufficiently wide class of multivariate function germs   with  similar
properties?

For a long time, I thought that the answer to this
question was negative. In this section the class of
$\Cal{SC}$-germs is defined. This provides an affirmative answer to
this question. \smallskip

A subset $A$ in a connected $k$-dimensional analytic manifold
$Y$ is called {\it meager} if there exists a countable set of
open domains $U_i\subset M$ and a countable collection of
proper analytic subsets $A_i\subset U_i$ in these domains such
that $A\subset\bigcup A_i$.

The following definition plays a key role in what follows.

\definition { Definition} A germ $f_a$ of an analytic function
at a point $a\in \Bbb C^n$ is an {\it
$\Cal {SC}$-germ} if the following condition is fulfilled.
For every connected complex analytic manifold $Y$, every
analytic map $G: Y\to \Bbb C^n$ and every preimage $b$ of
the point $a$, $G(b)=a$, there exists a meager set $A\subset Y$
such that, for every path $\gamma: [0,1]\to Y$ originating at
the point $b$, $\gamma(0)=b$ and intersecting the set $A$ at
most at the initial moment, $\gamma (t)\not\in A$ for $t>0$,
the germ $f_a$ admits an analytic continuation along the path
$G\circ\gamma: [0,1]\to \Bbb C^n$.
\enddefinition
The following lemma is obvious.
\proclaim {Lemma 1}
The class of $\Cal SC$-germs contains all germs of analytic functions on $\Bbb C^N\setminus \Sigma$ where $\Sigma$ is an analytic subset in $\Bbb C^N$ where $N$ is a natural number. In particular the class contains all analytic germs of
$\Cal S$-functions of one variable and all germs of
meromorphic functions of many variables.
\endproclaim

The proof of the following Theorem 2 uses the results on extendability of
multivalued analytic functions along their singular point sets
(see [1]).

\proclaim {Theorem 2 (on stability of the class of $\Cal {SC}$-germs)}
The class of $\Cal {SC}$-germs on $\Bbb C^n$ is stable under the
operation of taking the composition with $\Cal {SC}$-germs of
$m$-variable functions,
the operation of differentiation and integration. It is stable under
solving algebraic equations whose coefficients are $\Cal {SC}$-germs and under solving holonomic systems of linear  differential equations whose coefficients are $\Cal {SC}$-germs.

\endproclaim

Theorem 2 implies the following corollary.
\proclaim {Corollary 3} If a germ $f$ is not an $\Cal {SC}$-germ  then $f$ is strongly non representable by generalized quadratures.  In particular it cannot be a germ of a function belonging to a certain Liouvillian class.
\endproclaim

\subhead {8. Monodromy group of a $\Cal {SC}$-germs}
\endsubhead
The {\it monodromy group and the  monodromy pair} of a $\Cal {SC}$-germ $f_a$ can by defined in same way as for $\Cal S$-functions of one variable. By definition the set $\Sigma\subset \Bbb C^n$ of singular points of  $f_a$ is  meager set. Take any point $x_0\in \Bbb C^n\setminus \Sigma$ and consider the action of the fundamental group $\pi_1 (\Bbb C^n\setminus \Sigma,x_0)$ on the set $F_{x_0}$ of all germs equivalent to  the germ $f_a$. The {\it monodromy group} of $f_a$ is the image of the fundamental group under this action. The {\it monodromy pair} of $f_a$ is the pair $[\Gamma,\Gamma_0]_0$ where $\Gamma$ is the monodromy group and $\Gamma_0$ is the stationary subgroup of a germ  $f\in F_{x_0}$. Up to an isomorphism the monodromy group and the monodromy pair are independent of a choice  of the point $x_0$ and the germ $f$.

\remark {Remark} If a  $\Cal {SC}$-germ $f_a$ is defined at a singular point  $a\in \Sigma$  then  the {\it monodromy group of $f_a$ along  $\Sigma$ } is  defined: one can  consider continuations of $f_a$  along curves $\gamma$ belonging to $\Sigma$ and define a singular set $\Sigma_1\subset \Sigma$ for $f_a$ along $\Sigma$.  The monodromy group  of $f_a$ along $\Sigma$  corresponds to the action the fundamental group of $\pi_1(\Sigma\setminus \Sigma_1,x_1)$  on the set of germs at $x_1\in \Sigma \setminus \Sigma_1$ obtained by continuation of $f_a$ along $\Sigma$. If the point $a$ belongs to $\Sigma_1$ then one can define also a monodromy group of $f_a$ along $\Sigma_1$ and so on. Thus in the multidimensional case one can associate to an $\Cal {SC}$-germ an an hierarchy of monodromy groups. All these monodromy groups (and corresponding monodromy pairs)  appear in multidimensional topological Galois theory. But the monodromy group and the monodromy pair we discuss above are  most important for our purposes.
\endremark

\subhead {9. Stability of certain classes of $\Cal{SC}$-germs}\endsubhead
One can prove the following  theorems.
\proclaim {Theorem 4 (see [1])}   The class of all
$\Cal {SC}$-germs, having a solvable monodromy   is stable
under composition, arithmetic operations,integration and differentiation.This class contains all germs of basic elementary functions and all  germs of single valued  functions whose set of singular points is a proper analytic set.
\endproclaim

\proclaim {Theorem 5 (see [1])}  The class of all
$\Cal {SC}$-germs, having a $k$-solvable monodromy pair {\rm (see [1], [2])}  is stable
under composition, arithmetic operations, integration, differentiation and solution of algebraic equations of degree at most $k$. This class contains all germs of basic elementary functions and all  germs of single valued  functions whose set of singular points is a proper analytic set.
\endproclaim

\proclaim {Theorem 6 (see [1])}   The class of all
$\Cal {SC}$-germs, having an almost solvable monodromy pair {\rm (see [1], [2])}  is stable
under composition, arithmetic operations, integration, differentiation and solution of algebraic equations. This class contains all germs of basic elementary functions and all  germs of single valued  functions whose set of singular points is a proper analytic set.
\endproclaim

Theorems 4--6 imply  the  following corollaries.

\proclaim {Result on  quadratures} If the monodromy group of a $\Cal {SC}$-germ $f$ is not solvable, then $f$ is   strongly non representable by quadratures.
\endproclaim

\proclaim {Result on  $k$-quadratures} If the monodromy pair of a $\Cal {SC}$-germ $f$ is not $k$-solvable, then $f$ is   strongly non representable by $k$-quadratures.
\endproclaim

\proclaim {Result on  generalized quadratures} If the monodromy pair  of a $\Cal {SC}$-germ $f$ is not  almost solvable, then $f$ is   strongly non representable by generalized quadratures.
\endproclaim

\subhead {10. Solvability and non solvability of algebraic equation}\endsubhead
Consider an irreducible  algebraic equation
$$ P_n y^n+P_{n-1}y^{n-1}+\dots +P_0=0 \tag 2$$ whose coefficients  $P_n,\dots, P_0$ are polynomials  of $N$ complex variables $x_1,\dots, x_N$. Let $\Sigma\subset \Bbb C^N$ be the singular set of the equation (2)
defined by the equation $P_nJ=0$ where $J$ is the discriminant of the polynomial (2).

\proclaim{Theorem 7 (see [1], [2], [4])} Let $y_{x_0}$ be a germ of analytic function at a point $x_0\in \Bbb C^N\setminus \Sigma$ satisfying the equation {\rm (2)}. If the monogromy group of the equation {\rm (2)} is solvable (is $k$-solvable) then the germ $y_{x_0})$ is representable by radicals (is representable by $k$-radicals).
\endproclaim
According to Camille Jordan's theorem (see [4])  the  Galois group of the equation (2) over the field $\Cal R$ of rational functions of $x_1,\dots,x_N$ it is isomorphic to the  monodromy group of  this equation (2). Thus Theorem 7 follows from Galois theory (see [1], [4]).

\proclaim{Theorem 8 (see [1])} Let $y_{x_0}$ be a germ of analytic function at a point $x_0\in \Bbb C^N$ satisfying the  equation {\rm (2)}. If the monogromy group of the equation is not solvable (is not $k$-solvable) then the germ $y_{x_0}$ is strongly non representable by quadratures (is strongly non representable by $k$-quadratures).
\endproclaim
Theorem 8 follows from the results on quadratures and on $k$-quadratures from the previous section.

Consider the universal degree $n$ algebraic function $y(a_n,\dots, a_0)$ defined by the equation
$$a_ny^n+\dots+a_0=0. \tag 3$$

It is easy to see that the monodromy group of the equation (3) is isomorphic to the group $S_n$ of all permutations of $n$ element. For $n\geq 5$ the group $S_n$ is unsolvable and it is not $k$-solvable group for $k<n$. Thus Theorem 8 implies the following strongest known version of the Abel-–Ruffini Theorem.

\proclaim{Theorem 9 (a version of the Abel--Ruffini Theorem)} Let $y_{a}$ be a germ of analytic function at a point $a$ satisfying the universal degree $n\geq 5$ algebraic equation. If $n\geq 5$ then the germ $y_{a}$ is strongly non representable by $(n-1)$ quadratures. In particular the germ $y_{a}$ is strongly non representable by quadratures.
\endproclaim

\subhead {10. Solvability and non solvability of  holonomic systems of linear differential equations}\endsubhead  Consider a system of $N$ linear
differential equations $L_j(y)=0$, $j=1,\dots, N$,
$$
L_j(y)=\sum a_{i_1,\dots, i_n}\frac{\partial^{i_1+\dots
+i_n}y}{\partial x_1^{i_1} \dots
\partial x^{i_n}_n} =0, \tag 4
$$
on an unknown function $y$, whose coefficients $a_{i_1,\dots, i_n}$
are rational functions of $n$ complex variables $x_1$, $\dots$, $x_n$.
Assume that the system (4) is holonomic in $\Bbb C^n\setminus \Sigma_1$ where $\Sigma_1$ is the union of poles of the coefficients $a_{i_1,\dots, i_n}$. Let $\Sigma_2\subset \Bbb C^n\setminus \Sigma_1$ be the  singular hypersurface of a holonomic system (4).

Every germ $y_a$ of a solution of the system at a point $a\in \Bbb C^n\Sigma$ where $\Sigma=\Sigma_1\cup \Sigma_2$ admits an analytic continuation along every path
avoiding the hypersurface $\Sigma$ so the monodromy group of the system (4) is well-defined.

\proclaim {Theorem 10 (see [1])}
If the monodromy group of the holonomic system (4) is not solvable (not $k$-solvable, not almost solvable), then a germ $y_a$ of almost every solution at a point $a\in \Bbb C^n\setminus \Sigma$ is strongly non representable by quadratures (is strongly non representable by $k$-quadratures, is strongly non representable by generalized quadratures).

\endproclaim
Theorem 10 follows from the results on quadratures, on $k$-quadratures and on generalized quadratures from  section 9.

A holonomic system is said to be
{\it regular}, if near the singular set $\Sigma$
and near infinity the solutions of the system grow at most polynomially.

\medskip

\proclaim {Theorem 11 (see [1])} If the monodromy group of a regular holonomic system  is  solvable (is $k$-solvable, is almost solvable), then a germ $y_a$ of almost every solution at a point $a\in \Bbb C^n\setminus \Sigma$ is  representable by quadratures (is  representable by $k$-quadratures, is  representable by generalized quadratures).
\endproclaim

\subhead {11. Completely integrablesystems of linear differential equations with small coefficients}\endsubhead Consider a completely integrable system of linear differential equations of the
following form
 $$
 dy=Ay,\eqno{(5)}
 $$
where $y=y_1$, $\dots$, $y_N$ is an unknown vector-function, and $A$ is a
$(N\times N)$-matrix consisting of differential one-forms with rational coefficients
on the space $\Bbb C^n$ satisfying the condition of complete integrability
$dA+A\wedge A=0$ and having the following form:
 $$
 A=\sum _{i=1}^kA_i\frac{dl_i}{l_i},
 $$
where $A_i$ are constant matrices, and $l_i$ are linear (not necessarily homogeneous)
functions on $\Bbb C^n$.

If the matrices $A_i$ can be simultaneously reduced to the
triangular form, then system (5), as any completely integrable
triangular system, is solvable by quadratures. Of course, there
exist solvable systems that are not triangular. However, if the
matrices $A_i$ are sufficiently small, then there are no such
systems. Namely, the following theorem holds.

\proclaim {Theorem 12 (see [1])}
A system {\rm (5)} that does not reduce to the triangular form and
such that the matrices $A_i$ have sufficiently small norms is
unsolvable by generalized quadratures in the following strong sense. At every point $a\in \Bbb C^n$ where the matrix $A$ is regular, and for almost any germ $y_a=(y_1,\dots,y_N)_a$ of a vector-function satisfying the system {\rm (5)}, there is a component $(y_i)_a$ which is strongly non representable by generalized quadratures.
\endproclaim

Multidimensional Theorem 12 is similar to the one dimensional Corollary 20 from [2]. Their proofs (see [1]) are also similar.  We only need to replace the reference to the (one-dimensional)
Lappo-Danilevsky theory with the reference to the
multidimensional version of it from [5].

\subhead{7. Acknowledgement}\endsubhead
I would like to thank Michael Singer who invited me to write comments for a new edition of the classical J.F.~Ritt's book ``Integration in finite terms" [6]. This preprint was written as a part of these comments. I also am grateful to Fedor Kogan who edited my English.
\bigskip

\centerline {REFERENCES}
\medskip
[1]  A. Khovanskii, Topological Galois theory. Solvability and unsolvability of equations in
finite terms. Translated by Valentina Kiritchenko and Vladlen Timorin. Series: Springer
Monographs in Mathematics. Springer Berlin Heidelberg. 2014, XVIII, 305 pp. 6 illus.

[2] A. Khovanskii, One dimensional topological Galois theory. 2019.\linebreak
arXiv:1904.03341 [math.AG].

[3] M. van der Put, M. Singer, Galois theory of linear differential equations (Springer, Berlin/New York, 2003).

[4] A. Khovanskii, On representability of algebraic functions by radicals.
2019. arXiv:1903.08632 [math.AG].

[5] V.P. Leksin, Riemann--Hilbert problem for analytic families of representations. Math. Notes 50, 832--838 (1991).

[6] J. Ritt, Integration in finite terms. Liouville's theory of elementary methods, N.Y. Columbia Univ. Press. 1948.

\end